\documentclass[a4paper,11pt,twoside,english]{article}
\usepackage{babel}
\usepackage{latexsym,amsfonts}
\usepackage{amsthm}
\usepackage{amsmath}
\usepackage{mathtools}
\usepackage{paralist}

\usepackage{lastpage}
\usepackage{fancyhdr}
\fancypagestyle{plain}{

  \fancyhead[L,C]{}
  \fancyhead[R]{\today}
  \fancyfoot[R,L]{}
  \fancyfoot[C]{\thepage \ of \pageref{LastPage}}
  }

\title{Some remarks on generalised lush spaces}
\author{Jan-David Hardtke}
\date{}

\DeclareMathOperator{\lin}{span}
\DeclareMathOperator{\co}{co}
\DeclareMathOperator{\aco}{aco}

\providecommand{\ssq}{\subseteq}
\providecommand{\id}{\ensuremath{\mathrm{id}}}
\providecommand{\N}{\ensuremath{\mathbb{N}}}

\providecommand{\R}{\ensuremath{\mathbb{R}}}

\providecommand{\U}{\ensuremath{\mathcal{U}}}
\providecommand{\eps}{\ensuremath{\varepsilon}}
\providecommand{\cl}[2][]{\ensuremath{\overline{#2}^{#1}}}

\providecommand{\caco}[2][]{\ensuremath{\overline{\aco}^{#1} #2}}

\providecommand{\dist}[2]{\ensuremath{\operatorname{d} (#1,#2)}}

\providecommand{\keywords}[1]{
  {\let\thefootnote=\relax
  \footnote{{\em Keywords}: #1}}
  \addtocounter{footnote}{-1}
  }

\providecommand{\AMS}[1]{
  {\let\thefootnote=\relax
  \footnote{{\em AMS Subject Classification} (2010): #1}}
  \addtocounter{footnote}{-1}
  }

\providecommand{\address}{
  {\sc \noindent Department of Mathematics \\
  Freie Universit\"at Berlin \\
  Arnimallee 6, 14195 Berlin \\
  Germany \\}
  }

\DeclarePairedDelimiter{\set}{\lbrace}{\rbrace}
\DeclarePairedDelimiter{\paren}{\lparen}{\rparen}

\DeclarePairedDelimiter{\abs}{\lvert}{\rvert}
\DeclarePairedDelimiter{\norm}{\lVert}{\rVert}
\DeclarePairedDelimiter{\dotp}{\langle}{\rangle}

\theoremstyle{definition}
\newtheorem{definition}{Definition}[section]
\newtheorem*{definition*}{Definition}

\newtheorem*{example*}{Example}
\newtheorem{remark}[definition]{Remark}
\newtheorem*{remark*}{Remark}
\theoremstyle{plain}

\newtheorem*{lemma*}{Lemma}
\newtheorem{proposition}[definition]{Proposition}
\newtheorem*{proposition*}{Proposition}
\newtheorem{theorem}[definition]{Theorem}
\newtheorem*{theorem*}{Theorem}

\newtheorem*{corolary*}{Corollary}

\newenvironment{Proof}[1][\proofname]{\begin{proof}[#1] \setlength{\parindent}{0pt}}{\end{proof}}
\newenvironment{Abstract}{\centering\begin{minipage}{0.8\textwidth} \noindent \small {\sc Abstract.}}{\end{minipage}\par}

\usepackage{color}
\definecolor{darkgreen}{rgb}{0,0.5,0}

\numberwithin{equation}{section}
\newtagform{colored}[\color{blue}]{\color{blue}(}{\color{blue})}
\usetagform{colored}

\usepackage[colorlinks,linkcolor=blue,citecolor=red,urlcolor=darkgreen]{hyperref}
\providecommand{\email}{{\it E-mail address:} \href{mailto:hardtke@math.fu-berlin.de}{\tt hardtke@math.fu-berlin.de}}

\usepackage{amsrefs}

\begin{document}

\maketitle

\begin{Abstract}
X. Huang et al. recently introduced the notion of generalised lush (GL) spaces in \cite{huang}, which, 
at least for separable spaces, is a generalisation of the concept of lushness introduced in \cite{boyko1}. 
The main result of \cite{huang} is that every GL-space has the so called Mazur-Ulam property (MUP).\par
In this note, we will prove some properties of GL-spaces (further than those already established in \cite{huang}), 
for example, every $M$-ideal in a GL-space is again a GL-space, ultraproducts of GL-spaces are again GL-spaces, 
and if the bidual $X^{**}$ of a Banach space $X$ is GL, then $X$ itself still has the MUP.
\end{Abstract}
\AMS{46B20}

\section{Introduction}\label{sec:intro}
Our notation is as follows: if not otherwise stated, $X$ denotes a real Banach space, $X^*$ its dual, $B_X$ its closed 
unit ball and $S_X$ its unit sphere. For a subset $A$ of $X$, we denote by $\cl{A}$ its norm-closure and by $\co{A}$ resp. 
$\aco{A}$ its convex resp. absolutely convex hull. By $\dist{x}{A}$ we denote the distance from a point $x\in X$ to the 
set $A$. For any functional $x^*\in S_{X^*}$ and any $\eps>0$ let $S(x^*,\eps):=\set*{x\in B_X:x^*(x)>1-\eps}$ be the slice
of $B_X$ induced by $x^*$ and $\eps$.\par
Now let us begin by recalling the classical Mazur-Ulam-theorem (see \cite{mazur}), which states that every bijective isometry $T$ 
between two real normed spaces $X$ and $Y$ must be affine, i.\,e. $T(\lambda x+(1-\lambda)y)=\lambda T(x)+(1-\lambda)T(y)$
for all $x,y\in X$ and every $\lambda\in [0,1]$ (equivalently, $T-T(0)$ is linear). A simplified proof of this theorem was given
in \cite{vaeisaelae}. See also the recent paper \cite{nica} for an even further simplified argument.\par
In 1972, Mankiewicz \cite{mankiewicz} proved the following generalisation of the Mazur-Ulam theorem: if $A\ssq X$ and $B\ssq Y$ 
are convex with non-empty interior or open and connected, then every bijective isometry $T:A \rightarrow B$ can be extended to a 
bijective affine isometry $\tilde{T}:X \rightarrow Y$. This result implies in particular that every bijective isometry from $B_X$ 
onto $B_Y$ is the restriction of a linear isometry from $X$ onto $Y$. Tingley asked in \cite{tingley} whether the same is true if 
one replaces the unit balls of $X$ and $Y$ by their respective unit spheres. As a first step towards solving this problem, Tingley proved 
in \cite{tingley} that for finite-dimensional spaces $X$ and $Y$, every bijective isometry $T:S_X \rightarrow S_Y$ satisfies $T(-x)=-T(x)$ 
for all $x\in S_X$.\par
Though Tingley's problem remains open to the present day even in two dimensions, affirmative answers have been obtained for many special 
classes of spaces. In particular, the answer is ``yes'' if $Y$ is an (a priori) arbitrary Banach space and $X$ is any of the classical 
Banach spaces $\ell^p(I)$, $c_0(I)$, for $1\leq p \leq\infty$ and $I$ any index set, or $L^p(\mu)$, for $1\leq p\leq\infty$ and $\mu$ a 
$\sigma$-finite measure (see \cites{ding, fang2, ruidong, tan1, tan2} and further references therein). The answer is also known to be 
positive for $Y$ arbitrary and $X=C(K)$ if $K$ is a compact metric space (see \cite{fang1}).\par 
The notion of Mazur-Ulam property was introduced in \cite{cheng}: a (real) Banach space $X$ is said to have the Mazur-Ulam property (MUP) 
if for every Banach space $Y$ every bijective isometry between $S_X$ and $S_Y$ can be extended to a linear isometry between $X$ and $Y$.\par
Let us now recall that a Banach space $X$ is called a CL-space resp. an almost CL-space if for every maximally convex subset $F$ of $S_X$ one has $B_X=\aco{F}$ 
resp. $B_X=\caco{F}$. CL-spaces were introduced by Fullerton in \cite{fullerton}, almost CL-spaces were introduced by Lima (see \cites{lima1, lima2}).
Lima also proved that real $C(K)$ and $L^1(\mu)$ spaces (where $K$ is any compact Hausdorff space, $\mu$ any finite measure) are CL-spaces. The complex 
spaces $C(K)$ are also CL while $L^1(\mu)$ is in the complex case in general only almost CL (see \cite{martin}).\par
In \cite{cheng} Cheng and Dong proposed a proof that every CL-space whose unit sphere has a smooth point and every polyhedral space\footnote{A Banach 
space is called polyhedral if the unit ball of each of its finite-dimensional subspaces is a polyhedron, i.\,e. the convex hull of finitely many points.} 
has the MUP. Unfortunately this proof is not completely correct, as is mentioned in the introduction of \cite{kadets}. Kadets and Mart\'{\i}n proved 
in \cite{kadets} that every {\em finite-dimensional} polyhedral space has the MUP. In \cite{liu} Liu and Tan showed that every almost CL-space whose 
unit sphere admits a smooth point has the MUP.\par
Next let us recall the definition of lushness, which was introduced in \cite{boyko1} (in connection with a problem concerning
the numerical index of a Banach space). The space $X$ is said to be lush provided that for any two points $x,y\in S_X$ and 
every $\eps>0$ there exists a functional $x^*\in S_{X^*}$ such that $x\in S(x^*,\eps)$ and 
\begin{equation*}
\dist{y}{\aco{S(x^*,\eps)}}<\eps.
\end{equation*}
For example, every almost CL-space is lush but the converse is not true in general (see \cite{boyko1}*{Example 3.4}).\par
In \cite{huang} Huang, Liu and Tan proposed the following definition of generalised lush spaces: $X$ is called a generalised lush (GL) space if for every 
$x\in S_X$ and every $\eps>0$ there is some $x^*\in S_{X^*}$ such that $x\in S(x^*,\eps)$ and
\begin{equation*}
\dist{y}{S(x^*,\eps)}+\dist{y}{-S(x^*,\eps)}<2+\eps \ \ \forall y\in S_X.
\end{equation*}
It is proved in \cite{huang} that every almost CL-space and every separable lush space is a GL-space (see \cite{huang}*{Example 2.4} resp. 
\cite{huang}*{Example 2.5}). Also, in \cite{huang}*{Example 2.7} the space $\R^2$ equipped with the hexagonal norm $\norm{(x,y)}=\max\set*{\abs{y},\abs{x}+1/2\abs{y}}$
is given as an example of a GL-space which is not lush.\par
The following two propositions are proved in \cite{huang}.
\begin{proposition}[\cite{huang}*{Proposition 3.1}]\label{prop:huang1}
If $X$ is a GL-space, $Y$ any Banach space and $T:S_X \rightarrow S_Y$ is a (not necessarily onto) isometry, then
\begin{equation}\label{eq:1.1}
\norm{T(x)-\lambda T(y)}\geq\norm{x-\lambda y} \ \ \forall x,y\in S_X, \forall \lambda\geq 0.
\end{equation}
\end{proposition}

\begin{proposition}[\cite{huang}*{Proposition 3.3}]\label{prop:huang2}
If $X$ and $Y$ are Banach spaces and $T:S_X \rightarrow S_Y$ is an onto isometry which satisfies \eqref{eq:1.1}
then $T$ can be extended to a linear isometry from $X$ onto $Y$.
\end{proposition}

It follows that every GL-space (in particular, every almsot CL-space and every seperable lush space) has the MUP (\cite{huang}*{Theorem 3.2}).\par
The authors of \cite{huang} further call a Banach space $X$ a local GL-space if for every separable subspace $Y$ of $X$ there is a subspace $Z$
of $X$ which is GL and contains $Y$. Since lushness is separably determined (see \cite{boyko2}*{Theorem 4.2}) every lush space is a local GL-space 
(\cite{huang}*{Example 3.6}). From their Propositions 3.1 and 3.3 the authors of \cite{huang} conclude that even every local GL-space has the MUP 
(\cite{huang}*{Theorem 3.7}), thus {\em every} lush space (seperable or not) has the MUP (\cite{huang}*{Corollary 3.8}).\par
Many stability properties for GL-spaces have already been established in \cite{huang}, for example, if $X$ is GL then so is the space $C(K,X)$ of
all continuous functions from $K$ into $X$, where $K$ is any compact Huasdorff space (see \cite{huang}*{Theorem 2.10}). Also, the property GL is
preserved under $c_0$-, $\ell^1$- and $\ell^{\infty}$-sums (see \cite{huang}*{Theorem 2.11 and Proposition 2.12}). In the next section,
we will prove some further stability results.

\section{Stability results}\label{sec:stability}
\subsection{Ultraproducts}\label{subsec:ultra}
We begin with an easy observation on ultraproducts of GL-spaces. First we recall the definition of ultraproducts of Banach spaces (see for example \cite{heinrich}). 
Given a free ultrafilter $\U$ on $\N$, for every bounded sequence $(a_n)_{n\in \N}$ of real numbers there exists (by a compactness argument) a number 
$a\in \R$ such that for every $\eps>0$ one has $\set*{n\in \N:\abs{a_n-a}<\eps}\in\U$. Of course $a$ is uniquely determined. It is called the 
limit of $(a_n)_{n\in \N}$ along $\U$ and denoted by $\lim_{n, \U}a_n$.\par
Now for a given sequence $(X_n)_{n\in \N}$ of Banach spaces let us denote by $\ell^{\infty}((X_n)_{n\in \N})$ the space of all sequences
$(x_n)_{n\in \N}$ in the product $\prod_{n\in \N}X_n$ such that $\sup_{n\in \N}\norm{x_n}<\infty$. We put
\begin{align*}
&\mathcal{N}_{\U}:=\set*{(x_n)_{n\in \N}\in \ell^{\infty}((X_n)_{n\in \N}):\lim_{n, \U}\norm{x_n}=0} \ \mathrm{and} \\
&\prod_{n, \U}X_n:=\ell^{\infty}((X_n)_{n\in \N})/\mathcal{N}_{\U}.
\end{align*}
Equipped with the (well-defined) norm $\norm{[(x_n)_{n\in \N}]}_{\U}:=\lim_{n,\U}\norm{x_n}$ this quotient becomes a Banach space. It is 
called the ultraproduct of $(X_n)_{n\in \N}$ (with respect to $\U$). By the way, it is easy to see that the subspace $\mathcal{N}_{\U}$
is closed in $\ell^{\infty}((X_n)_{n\in \N})$ with respect to the usual sup-norm and that $\norm{\cdot}_{\U}$ coincides with the usual 
quotient-norm. For more information on ultraproducts the reader is referred to \cite{heinrich}.\par
In \cite{boyko2}*{Corollary 4.4} it is shown that the ultraproduct of a sequence of lush spaces is again lush, in fact it even satisfies a
stronger property, called ultra-lushness in \cite{boyko2}. We can easily prove an analogous result for GL-spaces. First we need a little remark.
\begin{remark}\label{rem:aux}
If $X$ is a GL-space, $x\in S_X$ and $\eps>0$ then there is some $x^*\in S_{X^*}$ such that $x\in S(x^*,\eps)$ and
\begin{equation*}
\dist{y}{S(x^*,\eps)}+\dist{y}{-S(x^*,\eps)}\leq (2+\eps)\norm{y}+2\abs*{1-\norm{y}} \ \ \forall y\in X.
\end{equation*}
\end{remark}

\begin{Proof}
Analogous to the proof of \cite{huang}*{Lemma 2.9}.
\end{Proof}

\begin{proposition}\label{prop:ultra}
Let $\U$ be a free ultrafilter on $\N$ and $(X_n)_{n\in \N}$ a sequence of GL-spaces. Let $Z=\prod_{n, \U}X_n$. Then the following holds: 
for every $z\in S_Z$ there is a functional $z^*\in S_{Z^*}$ with $z^*(z)=1$ such that for every $y\in S_Z$ there are $z_1,z_2\in S_Z$
with $z^*(z_1)=1=-z^*(z_2)$ and $\norm{y-z_1}+\norm{y-z_2}=2$. In particular, $Z$ is also a GL-space.
\end{proposition}

\begin{Proof}
Let $z=[(x_n)_{n\in \N}]\in S_Z$. Without loss of generality we may assume $x_n\neq 0$ for all $n\in \N$. By the previous remark we can find,
for every $n\in \N$, a functional $x_n^*\in S_{X_n^*}$ such that $x_n/\norm{x_n}\in S(x_n^*,2^{-n})$ and for every $v\in X_n$
\begin{equation}\label{eq:2.1}
\dist{v}{S(x_n^*,2^{-n})}+\dist{v}{-S(x_n^*,2^{-n})}\leq (2+2^{-n})\norm{v}+2\abs*{1-\norm{v}}.
\end{equation}
Define $z^*:Z \rightarrow \R$ by $z^*([(v_n)]):=\lim_{n,\U}x_n^*(v_n)$. Then $z^*$ is a well-defined element of $S_{Z^*}$ with $z^*(z)=1$
(because of $x_n^*(x_n)>(1-2^{-n})\norm{x_n}$ for all $n$).\par 
Now given any $y=[(y_n)]\in S_Z$ we can find, by \eqref{eq:2.1}, sequences $(u_n)_{n\in \N}$ and $(v_n)_{n\in \N}$ in $\ell^{\infty}((X_n)_{n\in \N})$
such that $-v_n,u_n\in S(x_n^*,2^{-n})$ and 
\begin{equation*}
\norm{u_n-y_n}+\norm{v_n-y_n}<(2+2^{-n})\norm{y_n}+2\abs{1-\norm{y_n}}+2^{-n}.
\end{equation*}
Also, because of $-v_n,u_n\in S(x_n^*,2^{-n})$, the sum on the left-hand side of the above equation is at least $2(1-2^{-n})$. Altogether it follows 
that $z_1:=[(u_n)]$ and $z_2:=[(v_n)]$ satisfy our requirements.
\end{Proof}

\subsection{F-Ideals}\label{subsec:F-ideals}
First we recall the following notions (see \cite{harmand}*{Chapter I, Definition 1.1}): a linear projection $P:X \rightarrow X$ is called an $M$-projection if
\begin{equation*}
\norm{x}=\max\set*{\norm{Px},\norm{x-Px}} \ \ \forall x\in X.
\end{equation*}
$P$ is called an $L$-projection if
\begin{equation*}
\norm{x}=\norm{Px}+\norm{x-Px} \ \ \forall x\in X.
\end{equation*}
A closed subspace $Y$ of $X$ is said to be an $M$-summand ($L$-summand) in $X$ if it is the range of some $M$-projection ($L$-projection) on $X$.
Equivalently, $Y$ is an $M$-summand ($L$-summand) in $X$ if and only if there is some closed subspace $Z$ in $X$ such that $X=Y\oplus_{\infty}Z$
($X=Y\oplus_1Z$). Also, $Y$ is called an $M$-ideal in $X$ if $Y^{\perp}$ is an $L$-summand in $X^*$ (where $Y^{\perp}:=\set*{x^*\in X^*:x^*|_Y=0}$
is the annihilator of $Y$).\par
Every $M$-summand is also an $M$-ideal, but not conversely. For example, if $K$ is a compact Hausdorff space and $A\ssq K$ is closed, then the subspace
$Y:=\set*{f\in C(K):f|_A=0}$ is always an $M$-ideal in $C(K)$ but it is an $M$-summand if and only if $A$ is also open in $K$ (see \cite{harmand}*{Chapter I, Example 1.4(a)}).\par
As is pointed out in \cite{harmand}, the notion of an ``$L$-ideal'' (i.\,e. a subspace whose annihilator is an $M$-summand in the dual) is not introduced because
every ``$L$-ideal'' is already an $L$-summand (see \cite{harmand}*{Chapter I, Theorem 1.9}).\par
Just to give a few more examples let us mention that $L^1(\mu)$ is an $L$-summand in its bidual for every $\sigma$-finite measure $\mu$ (cf. \cite{harmand}*{Chapter IV, Example 1.1(a)})
and, as can be found in \cite{harmand}*{Chapter III, Example 1.4(f)}, for a Hilbert space $H$ the space $K(H)$ of compact operators on $H$ is an $M$-ideal in $K(H)^{**}=L(H)$ 
(the space of all operators on $H$). For more information on $M$-ideals and $L$-summands the reader is referred to \cite{harmand}.\par
Of course it is also possible to consider more general types of summands and ideals (see the overview in \cite{harmand}*{p.45f} and the papers 
\cites{mena-jurado1, mena-jurado2, mena-jurado3, paya}, we will recall just the basic definitions here). Firstly, a norm $F$ on $\R^2$ is called
absolute if $F(a,b)=F(\abs{a},\abs{b})$ for all $(a,b)\in \R^2$ and it is called normalised if $F(1,0)=1=F(0,1)$. In the following, $F$ will always 
denote an absolute, normalised norm on $\R^2$. If $X$ and $Y$ are two Banach spaces, their $F$-sum $X\oplus_F Y$ is defined as the direct product 
$X\times Y$ equipped with the norm $\norm{(x,y)}=F(\norm{x},\norm{y})$, which is again a Banach space. For every $1\leq p\leq\infty$, the $p$-norm $F_p$ on
$\R^2$ is of course an absolute, normalised norm and the corresponding sum is just the usual $p$-sum of two Banach spaces.\par
An important  property of absolute, normalised norms is their monotoni\-city, i.\,e. for all $a,b,c,d\in \R$
\begin{equation*}
\abs{a}\leq\abs{c} \ \mathrm{and} \ \abs{b}\leq\abs{d} \ \Rightarrow \ F(a,b)\leq F(c,d).
\end{equation*}
A proof of this fact can be found for instance in \cite{bonsall}*{Lemma 2}. It follows in particular that $\abs{a},\abs{b}\leq F(a,b)$ holds for all $a,b\in \R$. 
We will use this later without further mention.\par
A linear projection $P:X \rightarrow X$ is called an $F$-projection if
\begin{equation*}
\norm{x}=F(\norm{Px},\norm{x-Px}) \ \ \forall x\in X
\end{equation*}
and of course, a closed subspace $Y$ of $X$ is said to be an $F$-summand in $X$ if it is the range of an $F$-projection 
(equivalently, $X=Y\oplus_F Z$ for some closed subspace $Z$). Finally, $Y$ is called an $F$-ideal if $Y^{\perp}$ is an $F^*$-summand in $X^*$,
where $F^*$ is the reversed dual norm of $F$, i.\,e.
\begin{equation*}
F^*(a,b)=\sup\set*{\abs{av+bu}:(u,v)\in \R^2 \ \mathrm{with} \ F(u,v)\leq 1} \ \ \forall (a,b)\in \R^2.
\end{equation*}
Then the $L$- resp. $M$-summands ($M$-ideals) are just the $F_1$- resp. $F_{\infty}$-summands ($F_{\infty}$-ideals). Every $F$-summand is also an $F$-ideal (see \cite{paya}*{Lemma 8}).
It is known that $F$-summands and $F$-ideals coincide (in every Banach space) if and only if the point $(0,1)$ is an extreme point of the unit ball of 
$(\R^2,F)$ (see \cite{paya}*{Corollary 10 and Remark 12} and the results in section 2 of \cite{mena-jurado1}).\par
It was proved in \cite{pipping} that every $L$-summand and every $M$-ideal in a lush space is again lush. In \cite{huang}*{Theorem 2.11} it is shown that the
$c_0$-sum of a family of Banach spaces is GL if {\em and only if} each summand is GL. So $M$-summands in GL-spaces are again GL. It is possible to extend this 
result to a class of $F$-ideals which includes in particular all $M$-ideals. The main tool of the proof is, as in \cite{pipping}, the principle of local reflexivity
(see \cite{albiac}*{Theorem 11.2.4}).
\begin{theorem}\label{thm:GLideal}
If $F$ is an absolute, normalised norm on $\R^2$ such that $(0,1)$ is an extreme point of the unit ball of $(\R^2,F^*)$, $X$ is a GL-space and $Y$ is
an $F$-ideal in $X$, then $Y$ is also a GL-space.
\end{theorem}

\begin{Proof}
Let $X^*=Y^{\perp}\oplus_{F^*} U$ for a suitable closed subspace $U\ssq X^*$. It easily follows that $U$ can be canonically identified with $X^*/Y^{\perp}$, 
which in turn can be canonically identified with $Y^*$, thus $X^*=Y^{\perp}\oplus_{F^*} Y^*$.\par
Now let $y\in S_Y$ and $0<\eps<1$ be arbitrary. Since $(0,1)$ is an extreme point of $B_{(\R^2,F^*)}$ and by an easy compactness argument there is a 
$0<\delta<\eps$ such that
\begin{equation}\label{eq:2.2}
F^*(a,b)=1 \ \mathrm{and} \ b\geq 1-\delta \ \Rightarrow \ \abs{a}\leq\eps.
\end{equation}
Because $X$ is GL we can find $x^*\in S_{X^*}$ such that $y\in S(x^*,\delta)$ and 
\begin{equation}\label{eq:2.3}
\dist{v}{S(x^*,\delta)}+\dist{v}{-S(x^*,\delta)}<2+\delta \ \ \forall v\in S_X.
\end{equation}
Write $x^*=(y^{\perp},y^*)$ with $y^{\perp}\in Y^{\perp}, y^*\in Y^*$ and $1=\norm{x^*}=F^*(\norm{y^{\perp}},\norm{y^*})$. Then $y^*(y)=x^*(y)>1-\delta>1-\eps$.
Since $\norm{y^*}\leq 1$ we get that $y\in S(y^*/\norm{y^*},\eps)$. It also follows that $\norm{y^*}>1-\delta$ and hence by \eqref{eq:2.2} we must have 
$\norm{y^{\perp}}\leq\eps$.\par
Next we fix an arbitrary $z\in S_Y$. By \eqref{eq:2.3} we can find $x_1\in S(x^*,\delta), x_2\in -S(x^*,\delta)$ such that
\begin{equation}\label{eq:2.4}
\norm{x_1-z}+\norm{x_2-z}<2+\delta.
\end{equation}
We have $X^{**}=Y^{**}\oplus_F (Y^{\perp})^*$, so if we consider $X$ canonically embedded 
in its bidual we can write $x_i=(y_i^{**},f_i)\in Y^{**}\oplus_F (Y^{\perp})^*$ for $i=1,2$. It follows that
\begin{equation*}
1-\delta<x^*(x_1)=f_1(y^{\perp})+y_1^{**}(y^*).
\end{equation*}
Taking into account that $\norm{y^{\perp}}\leq\eps$ and $\norm{f_1}\leq 1$ we obtain
\begin{equation}\label{eq:2.5}
y_1^{**}(y^*)>1-\delta-\eps>1-2\eps.
\end{equation}
Analogously one can see that
\begin{equation}\label{eq:2.6}
-y_2^{**}(y^*)>1-2\eps.
\end{equation}
It also follows from \eqref{eq:2.4} that
\begin{equation*}
F(\norm{y_1^{**}-z},\norm{f_1})+F(\norm{y_2^{**}-z},\norm{f_2})<2+\delta
\end{equation*}
and hence
\begin{equation}\label{eq:2.7}
\norm{y_1^{**}-z}+\norm{y_2^{**}-z}<2+\delta<2+\eps.
\end{equation}
We put $E=\lin\set*{y_1^{**},y_2^{**},z}$ and chose $0<\eta<2\eps$ such that
\begin{equation*}
\frac{1-2\eps}{1+\eta}>1-3\eps \ \ \mathrm{and} \ \ (1+\eta)(2+\eps)<2+2\eps.
\end{equation*}
Now the principle of local reflexivity (\cite{albiac}*{Theorem 11.2.4}) comes into play. It yields
a finite-dimensional subspace $V\ssq Y$ and an isomorphism $T:E \rightarrow V$ such that $\norm{T}, \norm{T^{-1}}\leq 1+\eta$, $T|_{E\cap Y}=\id$
and $y^*(Ty^{**})=y^{**}(y^*)$ for all $y^{**}\in E$. Let $y_i=Ty_i^{**}$ for $i=1,2$. Then $y^*(y_i)=y_i^{**}(y_i)$ and $\norm{y_i}\leq 1+\eta$.
By \eqref{eq:2.5}, \eqref{eq:2.6} and the choice of $\eta$ we obtain $\norm{y_i}>1-2\eps$ as well as
\begin{equation}\label{eq:2.8}
\frac{y_1}{\norm{y_1}}\in S\paren*{\frac{y^*}{\norm{y^*}},3\eps} \ \ \mathrm{and} \ \ \frac{y_2}{\norm{y_2}}\in -S\paren*{\frac{y^*}{\norm{y^*}},3\eps}.
\end{equation}
From \eqref{eq:2.7} and the choice of $\eta$ we get
\begin{equation*}
\norm{y_1-z}+\norm{y_2-z}=\norm{Ty_1^{**}-Tz}+\norm{Ty_2^{**}-Tz}<(1+\eta)(2+\eps)<2+2\eps.
\end{equation*}
Since $1-2\eps<\norm{y_i}\leq 1+\eta<1+2\eps$ we have $\norm{y_i-y_i/\norm{y_i}}<2\eps$ and thus it follows that
\begin{equation*}
\norm*{\frac{y_1}{\norm{y_1}}-z}+\norm*{\frac{y_2}{\norm{y_2}}-z}<2+6\eps,
\end{equation*}
which, in view of \eqref{eq:2.8}, finishes the proof.
\end{Proof}

As mentioned before, Theorem \ref{thm:GLideal} shows in particular that $M$-ideals in GL-spaces are again GL, which was proved for lushness in \cite{pipping}.
The proof of \cite{pipping} readily extends to the case of more general ideals that we considered above (we skip the details).
\begin{theorem}\label{thm:lushideal}
If $F$ is an absolute, normalised norm on $\R^2$ such that $(0,1)$ is an extreme point of the unit ball of $(\R^2,F^*)$, $X$ is a lush space and $Y$ is
an $F$-ideal in $X$, then $Y$ is also lush.
\end{theorem}

Theorem 2.11 in \cite{huang} also states that the $\ell^1$-sum of any family of Banach spaces is GL if {\em and only if} every summand is GL. 
The ``only if'' part of this statement just means that $L$-summands in GL-spaces are again GL-spaces. However, the proof of this part given
in \cite{huang} contains a slight mistake: the statement ``$\norm{u_{\lambda}}>1/2-\eps/2$ and $\norm{v_{\lambda}}>1/2-\eps/2$'' cannot be deduced 
from the two preceding lines (2.3) and (2.4) as claimed in \cite{huang}. For a counterexample just consider the sum $X:=\R\oplus_1\R$ and take 
$x:=(1,0)\in S_X$. Then the norm-one functional $x^*: X \rightarrow \R$ defined by $x^*(a,b):=a+b$ satisfies $x^*(x)=1$ and 
$\dist{y}{S}+\dist{y}{-S}=2$ for all $y\in S_X$, where $S:=\set*{z\in S_X:x^*(z)=1}$ (we even have $\aco{S}=B_X$). Now for $y:=(-1,0), u:=(u_1,u_2):=(0,1)$
and $v:=y$ we have $-v,u\in S$ and $\norm{y-u}_1+\norm{y-v}_1=2$. So if the claim in the proof of \cite{huang} was true we would obtain the
contradiction $\abs{u_1}\geq 1/2$.\par 
We will therefore include a slightly different proof for the inheritance of generalised lushness to $L$-summands here.
\begin{proposition}\label{prop:Lsum}
If $X$ is a GL-space and $Y$ is an $L$-summand in $X$, then $Y$ is also a GL-space.
\end{proposition}

\begin{Proof}
Write $X=Y\oplus_1Z$ for a suitable closed subspace $Z\ssq X$. Let $y\in S_Y$ and $0<\eps<1$. Take $0<\delta<\eps^2$. Since $X$ is GL there is a functional
$x^*=(y^*,z^*)$ in the unit sphere of $X^*=Y^*\oplus_{\infty}Z^*$ such that $y\in S(x^*,\delta)$ and
\begin{equation}\label{eq:2.9}
\dist{v}{S(x^*,\delta)}+\dist{v}{-S(x^*,\delta)}<2+\delta \ \ \forall v\in S_X.
\end{equation}
Since $x^*(y)=y^*(y)$ it follows that $y\in S(y^*/\norm{y^*},\delta)\ssq S(y^*/\norm{y^*},\eps)$.\par
Now fix an arbitrary $u\in S_Y$. Because of \eqref{eq:2.9} we can find $x_1\in S(x^*,\delta), x_2\in -S(x^*,\delta)$ such 
that $\norm{u-x_1}+\norm{u-x_2}<2+\delta$. Write $x_i=y_i+z_i$ with $y_i\in Y, z_i\in Z$ for $i=1,2$. It then follows that
\begin{equation}\label{eq:2.10}
\norm{u-y_1}+\norm{z_1}+\norm{u-y_2}+\norm{z_2}<2+\delta.
\end{equation}
We distinguish two cases. First we assume that $\norm{y_1}, \norm{y_2}\geq\eps$. Since $x_1\in S(x^*,\delta)$ we have that
\begin{equation*}
y^*(y_1)=x^*(x_1)-z^*(z_1)>1-\delta-\norm{z_1}\geq\norm{y_1}-\delta
\end{equation*}
and hence 
\begin{equation*}
y^*(\frac{y_1}{\norm{y_1}})>1-\frac{\delta}{\norm{y_1}}\geq1-\frac{\delta}{\eps}>1-\eps,
\end{equation*}
thus $y_1/\norm{y_1}\in S(y^*/\norm{y^*},\eps)$. Analogously one can see that $y_2/\norm{y_2}\in -S(y^*/\norm{y^*},\eps)$.
Furthermore, because of \eqref{eq:2.10} and since $\norm{y_i}+\norm{z_i}=\norm{x_i}>1-\delta$, we have
\begin{align*}
&\norm*{u-\frac{y_1}{\norm{y_1}}}+\norm*{u-\frac{y_2}{\norm{y_2}}}\leq\norm{u-y_1}+\norm{u-y_2}+\abs{1-\norm{y_1}}+\abs{1-\norm{y_2}} \\
&\leq\norm{u-y_1}+\norm{u-y_2}+\norm{z_1}+\norm{z_2}+2\delta<2+3\delta<2+3\eps.
\end{align*}
In the second case we have $\norm{y_1}<\eps$ or $\norm{y_2}<\eps$. If $\norm{y_1}<\eps$ it follows that $\norm{z_1}=\norm{x_1}-\norm{y_1}
>1-\delta-\eps>1-2\eps$ and hence, because of \eqref{eq:2.10},
\begin{equation*}
\norm{u-y_2}+\norm{z_2}<2+\delta-(1-2\eps)-\norm{u-y_1}<1+3\eps-(1-\norm{y_1})<4\eps.
\end{equation*}
Then in particular $\abs{y^*(u)-y^*(y_2)}<4\eps$ and thus (since $-x_2\in S(x^*\delta)$) we have
\begin{align*} 
&y^*(u)<4\eps+y^*(y_2)=4\eps+x^*(x_2)-z^*(z_2)<4\eps-(1-\delta)-z^*(z_2) \\
&<5\eps-1+\norm{z_2}\leq5\eps-\norm{y_2}\leq5\eps+\norm{u-y_2}-1<9\eps-1.
\end{align*}
Hence $-u\in S(y^*/\norm{y^*},9\eps)$. But then
\begin{equation*}
\dist{u}{S(y^*/\norm{y^*},9\eps)}+\dist{u}{-S(y^*/\norm{y^*},9\eps)}=\dist{u}{S(y^*/\norm{y^*},9\eps)}\leq 2.
\end{equation*}
If $\norm{y_2}<\eps$ an analogous argument shows that $u\in S(y^*/\norm{y^*},9\eps)$ and thus the proof is complete.
\end{Proof}

\subsection{Inheritance from the bidual}\label{subsec:inher bidual}
Next we would like to prove that every Banach space $X$ whose bidual is GL has itself the MUP (in fact we will prove a little bit more).
First consider the following (at least formal) weakening of the definition of GL-spaces.
\begin{definition}\label{def:property star}
A real Banach space $X$ is said to have the property $(*)$ provided that for every $\eps>0$ and all $x,y_1,y_2\in S_X$ there exists a
functional $x^*\in S_{X^*}$ such that $x\in S(x^*,\eps)$ and 
\begin{equation*}
\dist{y_i}{S(x^*,\eps)}+\dist{y_i}{-S(x^*,\eps)}<2+\eps \ \ \mathrm{for} \ i=1,2.
\end{equation*}
\end{definition}
The exact same proof as in \cite{huang} shows that \cite{huang}*{Proposition 3.1} (Proposi\-tion \ref{prop:huang1} in our introduction) holds
true not only for GL-spaces but for all spaces with property $(*)$ and consequently every space with property $(*)$ has the MUP.
Now we will prove that property $(*)$ inherits from $X^{**}$ to $X$ and thus in particular $X$ has the MUP if $X^{**}$ is a GL-space.
\begin{theorem}\label{thm:star bidual}
If $X^{**}$ has property $(*)$, then so does $X$.
\end{theorem}

\begin{Proof}
Again the princiüple of local reflexivity is the key to the proof. If we fix $x,y_1,y_2\in S_X$ and $\eps>0$ and consider $X$ canonically embedded
into its bidual, then since the latter has property $(*)$ we can find $x^{***}\in S_{X^{***}}$ such that $x\in S(x^{***},\eps)$ and 
$u_1^{**}, u_2^{**}\in S(x^{***},\eps), v_1^{**}, v_2^{**}\in -S(x^{***},\eps)$ with
\begin{equation}\label{eq:2.11}
\norm{y_i-u_i^{**}}+\norm{y_i-v_i^{**}}<2+\eps \ \ \mathrm{for} \ i=1,2.
\end{equation}
If we also consider $X^*$ canonically embedded into $X^{***}$ then by Goldstine's theorem $B_{X^*}$ is weak*-dense in $B_{X^{***}}$, so we can find
$\tilde{x}^*\in B_{X^*}$ such that
\begin{align*}
&\abs{u_i^{**}(\tilde{x}^*)-x^{***}(u_i^{**})}\leq\eps, \ \abs{v_i^{**}(\tilde{x}^*)-x^{***}(v_i^{**})}\leq\eps \ \ \mathrm{for} \ i=1,2 \\ 
&\mathrm{and} \ \ \abs{\tilde{x}^*(x)-x^{***}(x)}\leq\eps.
\end{align*}
We put $x^*=\tilde{x}^*/\norm{\tilde{x}^*}$. It follows that $x\in S(x^*,2\eps)$ as well as $x^*\in S(u_i^{**},2\eps)$ and $-x^*\in S(v_i^{**},2\eps)$ 
for $i=1,2$.\par
Now let $E:=\lin\set*{x,y_1,y_2,u_1^{**},u_2^{**},v_1^{**},v_2^{**}}\ssq X^{**}$ and choose $0<\delta<\eps$ such that
\begin{equation*}
\frac{1-2\eps}{1+\delta}>1-3\eps \ \ \mathrm{and} \ \ (2+\eps)(1+\delta)<2+2\eps.
\end{equation*}
By the principle of local reflexivity (\cite{albiac}*{Theorem 11.2.4}) there is a finite-dimensional subspace $F$ of $X$ and an isomorphism 
$T:E \rightarrow F$ such that $\norm{T}, \norm{T^{-1}}\leq 1+\delta$, $T|_{X\cap E}=\id$ and $x^*(Tx^{**})=x^{**}(x^*)$ for all $x^{**}\in E$.\par
Put $\tilde{u}_i:=Tu_i^{**}$ and $\tilde{v}_i:=Tv_i^{**}$ as well as $u_i=\tilde{u}_i/\norm{\tilde{u}_i}$ and $v_i=\tilde{v}_i/\norm{\tilde{v}_i}$ 
for $i=1,2$. We then have 
\begin{equation}\label{eq:2.12}
\frac{1-\eps}{1+\delta}\leq\norm{\tilde{u}_i}, \norm{\tilde{v}_i}\leq 1+\delta \ \ \mathrm{for} \ i=1,2.
\end{equation}
It follows that
\begin{equation*}
x^*(u_i)=\frac{x^*(Tu_i^{**})}{\norm{\tilde{u}_i}}\geq\frac{u_i^{**}(x^*)}{1+\delta}>\frac{1-2\eps}{1+\delta}>1-3\eps,
\end{equation*}
so $u_i\in S(x^*,3\eps)$ and similarly also $-v_i\in S(x^*,3\eps)$ for $i=1,2$. Furthermore, because of \eqref{eq:2.11}, we have
\begin{equation*}
\norm{y_i-\tilde{u}_i}+\norm{y_i-\tilde{v}_i}=\norm{Ty_i-Tu_i^{**}}+\norm{Ty_i-Tv_i^{**}}<(1+\delta)(2+\eps)<2+2\eps.
\end{equation*}
From \eqref{eq:2.12} we get that $\norm{u_i-\tilde{u}_i}, \norm{v_i-\tilde{v_i}}\leq\eps+\delta<2\eps$. Hence 
\begin{equation*}
\norm{y_i-u_i}+\norm{y_i-v_i}<2+6\eps \ \ \mathrm{for} \ i=1,2
\end{equation*}
and we are done.
\end{Proof}

By a similar argument one could also prove that lushness inherits from $X^{**}$ to $X$. This fact has already been established in
\cite{kadets0}*{Proposition 4.3}, albeit with a different proof (the proof in \cite{kadets0} is based on an equivalent formulation of 
lushness (\cite{kadets0}*{Proposition 2.1}) and does not use the principle of local reflexivity).

\section{GL-spaces and rotundity}
We start with an easy observation on Hilbert spaces.
\begin{remark}\label{rem:hilbert}
Let $H$ be a Hilbert space and put
\begin{equation*}
A:=\set*{(x,x^*)\in S_H\times S_{H^*}:x\in \ker x^*}.
\end{equation*}
Then 
\begin{equation}\label{eq:3.1}
\dist{x}{S(x^*,\eps)}=\sqrt{2(1-\sqrt{2\eps-\eps^2})}
\end{equation}
for all $0<\eps<1$ and all $(x,x^*)\in A$. Consequently,
\begin{equation*}
\lim_{\eps\to 0}\dist{x}{S(x^*,\eps)}=\sqrt{2} \ \ \text{uniformly in} \ (x,x^*)\in A.
\end{equation*}
\end{remark}

\begin{Proof}
Let $0<\eps<1$ and $(x,x^*)\in A$ be arbitrary. Take $y\in S_H$ such that $x^*=\dotp{\cdot,y}$. Put $U=\lin\set*{y}$ and $V=U^{\perp}=\ker x^*$.
Let $P_U$ and $P_V$ denote the orthogonal projections from $H$ onto $U$ and $V$ respectively. Since $x\in V$ we have $P_Vx=x$. It follows that 
for any $z\in S(x^*,\eps)$ we have
\begin{align*}
&\norm{x-z}^2=\norm{P_V(x-z)}^2+\norm{x-z-P_V(x-z)}^2 \\
&=\norm{x-P_Vz}^2+\norm{P_Vz-z}^2=\norm{x-P_V}^2+\norm{P_Uz}^2 \\
&=\norm{x-P_Vz}^2+\abs{\dotp{z,y}}^2=\norm{x-P_Vz}^2+\abs{x^*(z)}^2 \\
&\geq(1-\norm{P_Vz})^2+(1-\eps)^2.
\end{align*}
But $\norm{P_Vz}^2+\norm{P_Uz}^2=\norm{z}^2\leq 1$ and hence $\norm{P_Vz}^2\leq 1-\abs{x^*(z)}^2\leq 1-(1-\eps)^2$. Putting everything together we get
\begin{equation*}
\norm{x-z}^2\geq(1-\sqrt{1-(1-\eps)^2})^2+(1-\eps)^2=2(1-\sqrt{2\eps-\eps^2}),
\end{equation*}
which proves the ``$\geq$'' part of \eqref{eq:3.1}. On the other hand, if $0<\delta<\eps$ and we put $\lambda:=\sqrt{1-(1-\delta)^2}$ and 
$u:=\lambda x+(1-\delta)y$, then $\norm{u}^2=\lambda^2+(1-\delta)^2=1$ and $x^*(u)=\dotp{u,y}=1-\delta>1-\eps$. So $u\in S(x^*,\eps)$ and
thus $\dist{x}{S(x^*,\eps)}^2\leq\norm{x-u}^2=(1-\lambda)^2+(1-\delta)^2=2(1-\sqrt{2\delta-\delta^2})$. Letting $\delta\to \eps$ gives the 
desired result.
\end{Proof}

It immediately follows from Remark \ref{rem:hilbert} that a Hilbert space (with dimension at least two) is not a GL-space. It is possible to
generalise this statement in a certain sense. To do so let us first recall some basic rotundity notions (see for example \cite{fabian}*{Chapters 8--9}). 
A Banach space $X$ is said to be strictly convex if $x,y\in S_X$ and $\norm{x+y}=2$ already implies $x=y$. $X$ is called uniformly rotund (or uniformly 
convex) if for any two sequences $(x_n)_{n\in \N}$ and $(y_n)_{n\in \N}$ in the unit sphere of $X$ the condition $\norm{x_n+y_n}\to 2$ implies $\norm{x_n-y_n}\to 0$.
Finally, $X$ is said to be locally uniformly rotund (or locally uniformly convex) if for every point $x\in S_X$ one has that $x_n\in S_X$ for each 
$n\in \N$ and $\norm{x_n+x}\to 2$ implies $\norm{x_n-x}\to 0$. Such points $x$ are called LUR points of the $S_X$ (an easy normalisation argument 
shows that we can replace the condition $\norm{x_n}=1$ for all $n$ by $\norm{x_n}\to 1$ in the definition of LUR points).\par
Of course uniform rotundity implies local uniform rotundity which in turn implies strict convexity of the space, but the converses of these two 
implications are false in general (though by an easy compactness argument it can be seen that the three notions coincide for finite-dimensional Banach spaces). 
The most prominent examples of uniformly rotund spaces are the Hilbert spaces and, more generally, the spaces $L^p(\mu)$ for any measure $\mu$ and 
any $1<p<\infty$.\par
Now our aforementioned generalisation reads as follows.
\begin{proposition}\label{prop:LUR GL}
Let $X$ be an infinite-dimensional GL-pace with a (not necessarily countable) $1$-unconditional basis. Then $S_X$ has no LUR points.
\end{proposition}

\begin{Proof}
By fixing a normalised $1$-unconditional basis $(x_i)_{i\in I}$ of $X$ (where $I$ is a suitable index set) we can identify $X$ with a 
subpsace of $\R^I$ which contains all functions with finite-support and carries a norm $\norm{\cdot}$ such that $x\in X$ and $y\in \R^I$ with 
$\abs{y}\leq\abs{x}$ implies $y\in X$ with $\norm{y}\leq\norm{x}$, and $\norm{e_i}=1$ for all $i\in I$, where $e_i(j)=0$ for $i\neq j$ and $e_i(i)=1$. 
Note that $X$ is then contained in $c_0(I):=\set*{y\in \R^I:\set*{i\in I:\abs{y(i)}>\eps} \ \text{is finite for every} \ \eps>0}$.\par
Now let $x\in X$ be arbitrary. We wish to show that $x$ is not an LUR point of $S_X$. By replacing $x$ with $\abs{x}$ if necessary 
we may assume without loss of generality that $x(i)\geq 0$ for all $i\in I$.\par
If $x$ was an LUR point of $S_X$ then there would be $\eta>0$ such that for every $y\in X$
\begin{equation}\label{eq:3.2}
\abs*{\norm{y}-1}\leq\eta \ \mathrm{and} \ \abs*{\norm{x+y}-2}\leq\eta \ \Rightarrow \ \norm{x-y}\leq\frac{1}{2}.
\end{equation}
Let $\eps>0$ be arbitrary. Again since $x$ is assumed to be an LUR point there would be $0<\delta\leq\eps$ such that
\begin{equation*}
y\in B_X, \ \norm{x+y}\geq2(1-\delta) \ \Rightarrow \ \norm{x-y}\leq\eps.
\end{equation*}
Since $X$ is a GL-space there is some $x^*\in S_{X^*}$ with $x\in S(x^*,\delta)$ and
\begin{equation*}
\dist{y}{S(x^*,\delta)}+\dist{y}{-S(x^*,\delta)}<2+\delta \ \ \forall y\in S_X.
\end{equation*}
Choose $i_0\in I$ such that $x(i_0)\leq\eps$ and find $y_1\in S(x^*,\delta), y_2\in -S(x^*,\delta)$ with 
\begin{equation}\label{eq:3.3}
\norm{e_{i_0}-y_1}+\norm{e_{i_0}-y_2}<2+\delta.
\end{equation}
We have $\norm{x+y_1}\geq x^*(x+y_1)>2(1-\delta)$ and hence the choice of $\delta$ implies $\norm{x-y_1}\leq\eps$. Similarly, $\norm{x+y_2}\leq\eps$.
Combining this with \eqref{eq:3.3} gives
\begin{equation}\label{eq:3.4}
\norm{x+e_{i_0}}+\norm{x-e_{i_0}}<2+\delta+2\eps\leq2+3\eps.
\end{equation}
Note that $\abs*{e_{i_0}-x}=\abs*{x+e_{i_0}-2x(i_0)e_{i_0}}$, hence
\begin{equation*}
\abs*{\norm{e_{i_0}-x}-\norm{e_{i_0}+x}}\leq\norm{x+e_{i_0}-2x(i_0)-(x+e_{i_0})}=2x(i_0)\leq2\eps.
\end{equation*}
In view of \eqref{eq:3.4} it follows that
\begin{equation}\label{eq:3.5}
2\norm{x+e_{i_0}}<2+5\eps.
\end{equation}
Next we show that $\norm{x+e_{i_0}}>1+\eta$. If not, then since $\norm{x+e_{i_0}}\geq\norm{x}=1$ (by the monotonicity of $\norm{\cdot}$)
we would have $\abs*{\norm{x+e_{i_0}}-1}\leq\eta$ and also $2=2\norm{x}\leq\norm{2x+e_{i_0}}\leq2+\eta$, hence the choice of $\eta$ would imply 
$1=\norm{e_{i_0}}\leq1/2$. Thus $\norm{x+e_{i_0}}>1+\eta$.\par
But then by \eqref{eq:3.5} we get $\eta<5\eps/2$ and since $\eps$ was arbitrary it follows that $\eta\leq 0$. With this contradiction the proof is finished.
\end{Proof}

As regards finite-dimensional GL-spaces, virtually the same argument as in the previous proof can be used to show the following (in fact the argument 
is even simpler, for in finite-dimensional spaces the compactness of the unit ball allows us to take $\eta=\eps=\delta=0$ and $0$ instead $1/2$ at the 
end of \eqref{eq:3.2}).
\begin{proposition}\label{prop:finitedim GL}
Let $\norm{\cdot}$ be an absolute, normalised norm on $\R^n$, i.\,e. $\norm{(x_1,\dots,x_n)}=\norm{(\abs{x_1},\dots,\abs{x_n})}$ for all
$(x_1,\dots,x_n)\in \R^n$ and $\norm{e_i}=1$ for all $i=1,\dots,n$ (where $e_i$ is the canonical $i$-th basis vector). Let
\begin{equation*} 
A:=\set*{x=(x_1,\dots,x_n):\norm{x}=1 \ \mathrm{and} \ x_i=0 \ \mathrm{for\ some} \ i\in \set*{1,\dots,n}}.
\end{equation*}
If $(\R^n,\norm{\cdot})$ is a GL-space and $x\in A$, then $x$ is not an LUR point of the unit sphere of $(\R^n,\norm{\cdot})$.
\end{proposition}

\address
\email

\end{document}